\documentclass[a4paper]{amsart}
\usepackage{graphicx}

\newtheorem{thm}{Theorem}

\title{Minimal Grid Diagrams of the Prime Alternating Knots with 12 Crossings}

\author{Gyo Taek Jin}
\address{Department of Mathematical Sciences,
Korea Advanced Institute of Science and Technology,
Daejeon, 34141, Korea}
\email{trefoil@kaist.ac.kr}

\author{Hwa Jeong Lee}
\address{Department of Mathematics Education, Dongguk University-Gyeongju, Gyeongju 38066, Korea}
\email{hjwith@dongguk.ac.kr}

\def\knum#1#2#3#4#5{\put(10,145){\small #1#2#3}}

\begin{document}

\begin{abstract}
In this article, we give a list of minimal grid diagrams of the 12 crossing prime alternating knots. This is a continuation of the work in~\cite{Jin2020}.
\end{abstract}

\maketitle

\section{Introduction}

\begin{figure}[h!]
\bigskip\setlength{\unitlength}{0.5pt}
\begin{picture}(140,160)(0,0)
\put(10,70){\line(0,1){40}}
\put(80,10){\line(1,0){50}}
\put(20,60){\line(0,1){20}}
\put(110,20){\line(1,0){30}}
\put(30,100){\line(0,1){30}}
\put(100,30){\line(1,0){20}}
\put(40,90){\line(0,1){20}}
\put(50,40){\line(1,0){40}}
\put(50,40){\line(0,1){60}}
\put(110,50){\line(1,0){20}}
\put(60,120){\line(0,1){20}}
\put(20,60){\line(1,0){100}}
\put(70,70){\line(0,1){60}}
\put(10,70){\line(1,0){60}}
\put(80,10){\line(0,1){110}}
\put(20,80){\line(1,0){120}}
\put(90,40){\line(0,1){100}}
\put(40,90){\line(1,0){60}}
\put(100,30){\line(0,1){60}}
\put(30,100){\line(1,0){20}}
\put(110,20){\line(0,1){30}}
\put(10,110){\line(1,0){30}}
\put(120,30){\line(0,1){30}}
\put(60,120){\line(1,0){20}}
\put(130,10){\line(0,1){40}}
\put(30,130){\line(1,0){40}}
\put(140,20){\line(0,1){60}}
\put(60,140){\line(1,0){30}}
\end{picture}
\caption{A minimal grid diagram of the knot $12a1$}\label{fig:12a1grid}
\end{figure}
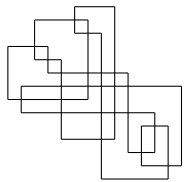

There are 1288 prime alternating knots with 12 crossings~\cite{knotinfo}.
We wrote a series of computer programs on Maple to find their minimal grid diagrams.
A \emph{grid diagram} is a knot diagram with finitely many horizontal segments and the same number of vertical segments such that the vertical segments cross over the horizontal segments at all crossings. See Figure\,\ref{fig:12a1grid}.

Imagine there is a vertical axis behind a grid diagram of a knot $K$ with $n$ vertical segments and $n$ horizontal segments. For each horizontal segment take the point on the axis at the same level to form a triangle and then replace the horizontal segment with the other two sides of the triangle. This results a polygonal embedding of $K$ with $3n$ edges such that each of  the $n$ half planes determined by the $n$ vertical segments and the axis 
contains a single arc of the knot. Namely, we have an \emph{arc presentation} of $K$ in $n$ arcs~\cite{C1995}. The minimal number of vertical segments in all grid diagrams of a knot is called the \emph{arc index}.
According to Theorem\,\ref{thm:Bae-Park}, the minimal grid diagrams of the 12 crossing knots have  14 vertical segments.
\begin{thm}[Bae-Park \cite{BP2000}]\label{thm:Bae-Park}
If $K$ is an alternating knot, then the arc index of $K$ is the minimal crossing number plus two.
\end{thm}

 For each of the mentioned 1288  knots, we pursued the following steps:
\begin{enumerate}
\item Using the DT notation from  \cite{knotinfo}, the regions divided by the knot diagram are described.
\item Considering the knot diagram as a planar graph, a spanning tree is generated whose contraction leads to an arc presentation of the knot.
\item The arc presentation  is converted to a grid diagram in various forms such as 2d and 3d graphics, latex picture commands and a sequence of 3d coordinates of vertices as a polygonal knot.
\item Using Knotplot~\cite{knotplot}, a DT notation for the grid diagram is obtained.
\item Using Knotscape~\cite{knotscape}, the grid diagram is confirmed to be the same as the original knot up to taking mirror images.
\end{enumerate}

\section{An example}
\begin{figure}[b!]
\includegraphics[width=0.55\textwidth]{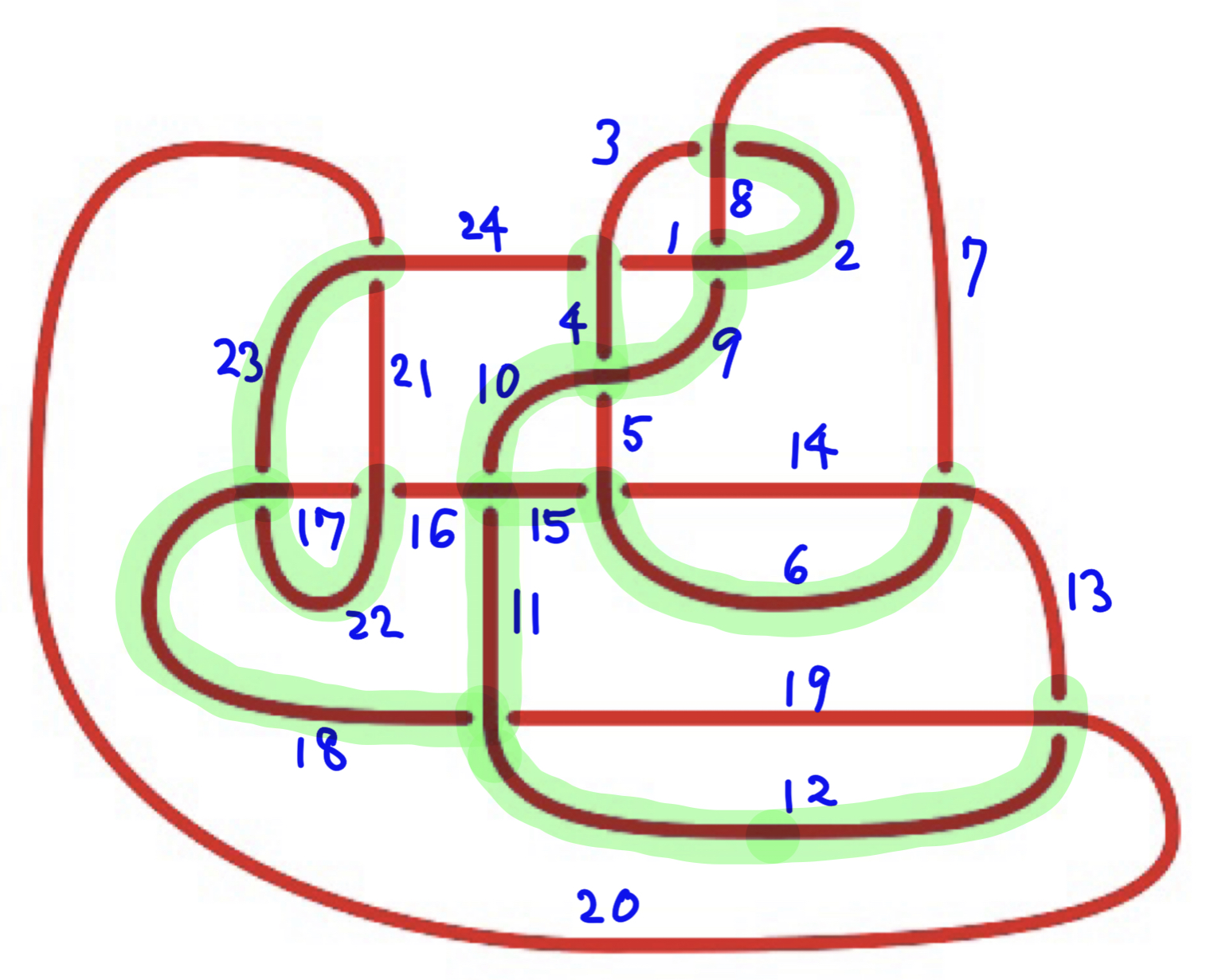}%
\includegraphics[width=0.55\textwidth]{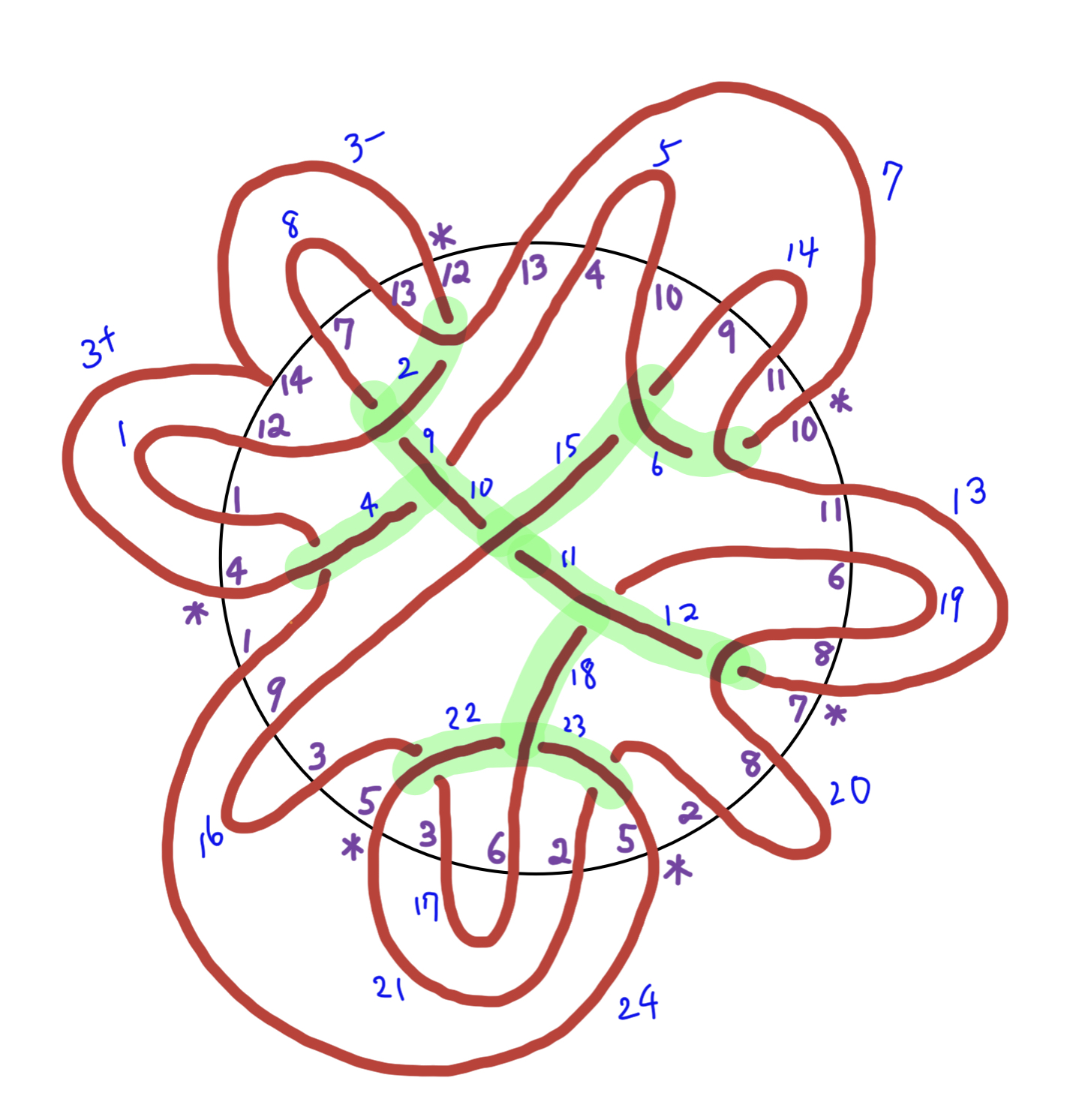}

\hfill(a)\hfill\hfill (b)\hfill\hfill
\caption{$12a1$}\label{fig:12a1}
\end{figure}
We give a  detailed description of the steps (1)--(3) using the knot $12a1$.
The diagram in Figure\,\ref{fig:12a1}(a) shows a minimal diagram of the knot $12a1$ obtained from \cite{knotinfo} with edge labels
compitible with the DT notation 	[4, 8, 10, 14, 2, 16, 20, 6, 22, 12, 24, 18] in the following sense. 
At a crossing incident to the four edges with lables $i$, $i+1$, $j$, $j+1$, modulo $24$, we may assume that $i$ and $j+1$ are odd and $i+1$ and $j$ are even. 
Then $i+1$ appears in the $(j/2+1)$-st (modulo 12) place in the DT notation. We orient the knot 
in the direction as the labels increase.
The diagram divides the plane into 14 regions, namely,
$[2, 8]$, $[6, 14]$, $[12, -19]$, $[17, -22]$, $[-21, -23, -17]$, $[1, -8, 3]$, $[4, -9, -1]$, $[10, -15, -5]$, $[-18, 23, -20, -12]$, $[22, 18, -11, 16]$, $[-24, 21, -16, -10, -4]$, $[9, 5, -14, 7, -2]$, $[15, 11, 19, 13, -6]$, and $[20, 24, -3, -7, -13]$ which are named by their oriented boundary edges which are compatible with their neighbors. For example, $[20, 24, -3, -7, -13]$ is the unbounded region whose boundary is oriented clockwise. The numbers with minus signs indicate that the corresponding edges are oriented in the opposite direction of the knot. All other regions are oriented counterclockwise.

The thickened edges of Figure\,\ref{fig:12a1}(a) form a spanning tree of the diagram. This tree is rooted at  the crossing between the edges 12 and 13. From the root, it grows as the following sequence of oriented edges  indicates.
\begin{equation}\label{eq:tree}
-12, -11, -18, -15, -10, -22, 23, 6, -9, -4, 2\tag{$\star$}
\end{equation}

The knot diagram in Figure\,\ref{fig:12a1}(b) is obtained from that of Figure\,\ref{fig:12a1}(a) by a plane isotopy. The edges are labeled in the same manner except that the one labeled with $3$ is divided into two parts. The thirteen arcs inside the circle are placed at distinct horizontal levels with heights as marked at their ends near the circle.
At the root of the tree the undercrossing arc has the height of 7 and the overcrossing arc 8. As one adds the edges of the tree in the order given by (\ref{eq:tree}),  new crossings appear one by one. At each new crossing, if a new arc undercrosses, then a lower height is given and if overcrosses then a higher height. 
The heights of the arcs in the circle decide the heights of the two ends of the thirteen edges outside the circle.
Suppose that the circle is the projection of the cylinder $x^2+y^2=1$ onto the $xy$-plane. If we collapse the cylinder to the $z$-axis, then the endpoints of the edges outside the circle have endpoints on the $z$-axis. Except the edge labeled 3, the other  twelve edges can be moved to be contained in distinct vertical half planes along the $z$-axis. We bring the center of the edge labeled 3 to the $z$-axis at a new height, highest as in the figure or lowest. Then the two parts can be placed in two new half planes.  The edge labeled 3 is handled in this special way because it is located at the extension of the last edge of the spanning tree. From the point on the circle between the edges 13 and 20, we read the heights of the endpoints in pairs counterclockwise. If the two ends on an edge are not adjacent, we read them at the $*$ marks. In this way we obtained
\begin{equation}\label{eq:12a1wheel}
\begin{aligned}
&[7, 11], [6, 8], [10, 13], [9, 11], [4, 10], [12, 14], [7, 13], \\
&[1, 12], [4, 14], [3, 9], [2, 5], [3, 6], [1, 5], [2, 8]
\end{aligned}\tag{$\star\star$}
\end{equation}
This sequence of pairs determines the grid diagram of Figure\,\ref{fig:12a1grid}. The first vertical segment spans the interval $[7, 11]$, the second spans $[6, 8]$, and so on.
Any cyclic permutation of (\ref{eq:12a1wheel}) gives another minimal grid diagram of the knot 12a1.

The idea behind this construction is originally due to Bae and Park~\cite{BP2000}, modified by the authors using the spanning tree~\cite{Jin2012}.
The spanning tree (\ref{eq:tree}) is constructed so that  the edges outside the circle have endpoints not interleaved among others except the one divided into two parts. 
Interested readers are encouraged to read \cite{BP2000} or \cite{Jin2012}.

\section{Minimal grid diagrams of 12 crossing prime alternating knots}
\vskip10pt
\setlength{\unitlength}{0.35pt}
\noindent


\section*{Acknowledgments}
The first author was supported by Basic Science Research Program through the National Research Foundation of Korea(NRF) funded by the Ministry of Education(2019R1A6A1A10073887). The corresponding author was supported by the National Research Foundation of Korea(NRF) grant funded by the Korea government (MSIT) (2019R1A2C1005506) and the Dongguk University Research Fund of 2020.

\end{document}